\documentclass{gtart_h}

\def\ifplaintex{\expandafter\ifx\csname documentclass\endcsname\relax}

\ifplaintex 
\else
\fi

\def\gt{{\mathsurround=0pt\it $\cal G\mskip-2mu$eometry \&\ 
$\cal T\!\!$opology}}        

\def\gtm{{\mathsurround=0pt\it $\cal G\mskip-2mu$eometry \&\ 
$\cal T\!\!$opology $\cal M\mskip-1mu$onographs}}    

\def\gtp{{\mathsurround=0pt\it $\cal G\mskip-2mu$eometry \&\ 
$\cal T\!\!$opology $\cal P\!$ublications}}  

\def\recd{{\small Received:\qua\receiveddate\ifx\reviseddate\relax
\else\qquad Revised:\qua\reviseddate\fi\par}} 

\def\volumenumber#1{\def\thevolumenumber{#1}}
\def\volumeyear#1{\def\thevolumeyear{#1}}
\def\volumename#1{\def\thevolumename{#1}}
\def\papernumber#1{\def\thepapernumber{#1}}
\def\pagenumbers#1#2{\def\startpage{#1}\def\finishpage{#2}}
\def\published#1{\def\publishdate{#1}}
\def\received#1{\def\receiveddate{#1}}
\def\revised#1{\def\reviseddate{#1}}
\def\accepted#1{\def\accepteddate{#1}}

\long\def\asciiabstract#1{\long\def\theasciiabstract{#1}}

\let\\\par
\let\thevolumenumber\relax\let\thepapernumber\relax
\let\thevolumeyear\relax\let\startpage\relax
\let\finishpage\relax\let\publishdate\relax\let\receiveddate\relax
\let\reviseddate\relax\let\accepteddate\relax\let\theasciititle\relax
\let\theasciiauthors\relax
\let\theasciiabstract\relax

\let\theerratum\relax\let\theasciiemail\relax
\let\theshortauthors\relax\let\theshorttitle\relax

\def\startpage{1}\def\finishpage{15}\def\thepapernumber{77}

\volumenumber{2}
\volumename{Proceedings of the Kirbyfest}
\volumeyear{1999}

\long\def\maketitlep{   

\count0=\startpage

\gtm\nl        
{\small Volume \thevolumenumber: \thevolumename\nl 
\ifx\theerratum\relax\else Erratum \erratumnumber\nl\fi
Pages \startpage--\finishpage\nl}

\vglue 0.1truein   

{\parskip=0pt\leftskip 0pt plus 1fil\def\\{\par\smallskip}{\ifplaintex\large
\else\Large\fi\bf\thetitle}\par\medskip}   
\vglue 0.05truein 

{\parskip=0pt\leftskip 0pt plus 1fil\def\\{\par}{\sc\theauthors}
\par\medskip}%
 
\vglue 0.03truein 

{\small\leftskip 25pt\rightskip 25pt{\bf Abstract}\stdspace\theabstract

{\bf AMS Classification}\stdspace\theprimaryclass
\ifx\thesecondaryclass\relax\else; \thesecondaryclass\fi\par
{\bf Keywords}\stdspace \thekeywords\par}\vglue 7pt

}   

\font\phead=cmsl9 scaled 950
\font=cmsl9 scaled 1050
\font\pnum=cmbx10 scaled 913
\font\lnum=cmbx10 
\font\pfoot=cmsl9 scaled 950
\font=cmsl9 scaled 1050
\ifplaintex
\headline{\vbox to 0pt{\vskip -4.5mm\line{\small\phead\ifnum
\count0=\startpage ISSN 1464-8997 (on line)
1464-8989 (printed) \hfill {\pnum\folio}\else\ifodd\count0\def\\{ }%
\ifx\theshorttitle\relax\thetitle\else\theshorttitle\fi\hfill{\pnum\folio}
\else\def\\{ and }{\pnum\folio}\hfill\ifx\theshortauthors\relax\theauthors
\else\theshortauthors\fi\fi\fi}\vss}}
\footline{\vbox to 0pt{\vglue 0mm\line{\small\pfoot\ifnum\count0=\startpage
Published \publishdate:\qua\copyright\ \gtp\hfill\else
\gtm, Volume \thevolumenumber\ (\thevolumeyear)\hfill\fi}\vss
}}
\else
\makeatletter
\def\@oddhead{{\small\lhead\ifnum\count0=\startpage ISSN 1464-8997 (on line)
1464-8989 (printed) \hfill {\lnum\number\count0}\else\ifodd\count0
\def\\{ }\ifx\theshorttitle\relax \thetitle \else\theshorttitle\fi\hfill
{\lnum\number\count0}\else\def\\{ and }{\lnum\number\count0}
\hfill\ifx\theshortauthors\relax 
\theauthors\else\theshortauthors\fi\fi\fi}}\def\@evenhead{@oddhead}
\def\@oddfoot{\small\lfoot\ifnum\count0=\startpage Published \publishdate:\qua\copyright\ \gtp\hfill\else
\gtm, Volume \thevolumenumber\ (\thevolumeyear)\hfill\fi}
\def\@evenfoot{@oddfoot}
\makeatother
\fi

\let\maketitlepage\maketitlep
\let\makeshorttitle\maketitlepage
\let\maketitle\maketitlepage

\newwrite\gtoutfile
\long\gdef\makeheadfile{  
{\def\\{, }\def\s{ }
\immediate\openout\gtoutfile head.xxx
\immediate\write\gtoutfile{Proxy-for: \ifx\theasciiauthors\relax
\theauthors\else\theasciiauthors\fi\s<\ifx\theasciiemail\relax\theemail\else\theasciiemail\fi>}
\immediate\write\gtoutfile{\noexpand\\}
\immediate\write\gtoutfile{Authors: \ifx\theasciiauthors\relax
\theauthors\else\theasciiauthors\fi}
{\def\\{ }\immediate\write\gtoutfile{Title: \ifx\theasciititle\relax
\thetitle\else\theasciititle\fi}}
\immediate\write\gtoutfile{Subj-class: GT or SG, GR etc}
\immediate\write\gtoutfile{MSC-class: \theprimaryclass\ifx\thesecondaryclass\relax\else, \thesecondaryclass\fi}
\immediate\write\gtoutfile{Journal-ref: Geom. Topol. Monogr. \thevolumenumber\s
(\thevolumeyear) \startpage-\finishpage}
\immediate\write\gtoutfile{Comments: Published by Geometry and Topology Monographs at}
\immediate\write\gtoutfile{\s\s\s  http://www.maths.warwick.ac.uk/gt/GTMon\thevolumenumber/paper\thepapernumber.abs.html}
\immediate\write\gtoutfile{\noexpand\\}
\immediate\write\gtoutfile{}
\ifx\theasciiabstract\relax
\immediate\write\gtoutfile{\theabstract}\else
\immediate\write\gtoutfile{\theasciiabstract}\fi
\immediate\write\gtoutfile{}
\immediate\write\gtoutfile{\noexpand\\}
\immediate\write\gtoutfile{}
\immediate\closeout\gtoutfile}}  

\def\maketitlepage{\maketitlep\makeheadfile}
\let\makeshorttitle\maketitlepage
\let\maketitle\maketitlepage

\volumenumber{4}
\volumename{Invariants of knots and 3-manifolds (Kyoto 2001)}
\volumeyear{2002}
\papernumber{22}
\pagenumbers{337}{362}
\received{25 November 2003}
\revised{24 January 2004}
\accepted{2 February 2004}
\published{3 February 2004}

\usepackage{amssymb,amsmath,latexsym}
\usepackage{cite,graphicx,calc,booktabs,picins,hyphenat}

\usepackage[dvips,bookmarks,bookmarksnumbered,pdfstartview=FitBH]{hyperref}

\def\mathcenter#1{%
  \vcenter{\hbox{#1}}%
}

\def\graph#1{
        \includegraphics[trim=-2 -2 -2 -2]{#1}
}

\def\mathgraph#1{
        \mathcenter{\graph{#1}}
}

\def\mathgraphsmall#1{
        \mathcenter{\includegraphics{#1}}
}

\newcommand{\RR}{\mathbb R}

\newcommand{\ZZ}{\mathbb Z}

\theoremstyle{plain}
\newtheorem{theorem}{Theorem}
\newtheorem{proposition}{Proposition}

\newtheorem{conjecture}[proposition]{Conjecture}
\newtheorem{observation}[proposition]{Observation}

\theoremstyle{definition}
\newtheorem{definition}[proposition]{Definition}
\newtheorem{exercise}[proposition]{Exercise}
\newtheorem{question}[proposition]{Question}

\theoremstyle{remark}

\newtheorem*{remark}{Remark}

\newcommand{\crosstet}{\smash{\mathgraphsmall{draws/misc.2}}}
\newcommand{\thetagraph}{\smash{\mathgraphsmall{draws/misc.1}}}
\newcommand{\tripletheta}{\mathgraphsmall{draws/misc.0}}

\begin{document}
\title{The algebra of knotted trivalent graphs\\and
  Turaev's shadow world}

\author{Dylan P. Thurston}
\address{Department of Mathematics, Harvard University\\Cambridge, 
MA 02138, USA}
\email{dpt@math.harvard.edu}

\begin{abstract}
  Knotted trivalent graphs (KTGs) form a rich algebra with a few
  simple operations: connected sum, unzip, and bubbling.  With these
  operations, KTGs are generated by the unknotted tetrahedron and
  M{\"o}bius strips.  Many previously known representations of knots,
  including knot diagrams and non-associative tangles, can be turned
  into KTG presentations in a natural way.
  
  Often two sequences of KTG operations produce the same output on all
  inputs.  These ``elementary'' relations can be subtle: for instance,
  there is a planar algebra of KTGs with a distinguished cycle.
  Studying these relations naturally leads us to Turaev's \emph{shadow
    surfaces}, a combinatorial representation of 3-manifolds based on
  simple 2-spines of 4-manifolds.  We consider the knotted trivalent
  graphs as the boundary of a such a simple spine of the 4-ball, and
  to consider a Morse-theoretic sweepout of the spine as a ``movie''
  of the knotted graph as it evolves according to the KTG operations.
  For every KTG presentation of a knot we can construct such a movie.
  Two sequences of KTG operations that yield the same surface are
  topologically equivalent, although the converse is not quite true.
\end{abstract}

\asciiabstract{%
  Knotted trivalent graphs (KTGs) form a rich algebra with a few
  simple operations: connected sum, unzip, and bubbling.  With these
  operations, KTGs are generated by the unknotted tetrahedron and
  Moebius strips.  Many previously known representations of knots,
  including knot diagrams and non-associative tangles, can be turned
  into KTG presentations in a natural way.
  
  Often two sequences of KTG operations produce the same output on all
  inputs.  These `elementary' relations can be subtle: for instance,
  there is a planar algebra of KTGs with a distinguished cycle.
  Studying these relations naturally leads us to Turaev's shadow
    surfaces, a combinatorial representation of 3-manifolds based on
  simple 2-spines of 4-manifolds.  We consider the knotted trivalent
  graphs as the boundary of a such a simple spine of the 4-ball, and
  to consider a Morse-theoretic sweepout of the spine as a `movie'
  of the knotted graph as it evolves according to the KTG operations.
  For every KTG presentation of a knot we can construct such a movie.
  Two sequences of KTG operations that yield the same surface are
  topologically equivalent, although the converse is not quite true.}

\maketitle

\section{Introduction}
\label{sec:intro}
In this paper we study the algebra of \emph{knotted trivalent graphs}
(KTGs).  A knotted trivalent graph is a framed\footnote{See
  Section~\ref{sec:KTG} for the precise notion of framing we use.}
embedding of a trivalent graph into $\RR^3$, modulo isotopy.  These
KTGs support some simple operations, forming an algebra-like structure.
Every knot may be presented as a sequence of
KTG operations starting with elementary graphs.  Thus we may use KTGs
as a novel representation of knot theory via generators and relations.

We may compare a KTG presentation with other representations of knots,
such as:
\begin{itemize}
\item Planar knot diagrams;
\item Braid closures;
\item $n$-bridge representations;
\item Pretzel representations;
\item Rational tangles and algebraic knots~\cite{Conway:Enumeration};
\item Parenthesized tangles~\cite{BarNatan:NAT,LM:FramedOriented}; and
\item Curves in a book with three pages%
\footnote{A book with three pages looks like this: $\mathgraph{draws/misc.30}$}~\cite{Dynnikov:3Page}.
\end{itemize}
These representations of knots all take an inherently 3-dimensional
object (a knot) and squash it into 2 dimensions (as in a knot diagram)
or sometimes even into 1 dimension (as in a parenthesized tangle).
The algebra of KTGs deals more directly with knots as
3-dimensional objects, while strictly generalizing all of the above
representations of knots, in the sense that, for instance, any knot
diagram can be turned into a sequence of operations on KTGs in a
natural way.\footnote{Here ``algebra'' is used in the universal
  algebra sense of a set supporting some operations}
We will illustrate how KTGs generalize knot diagrams in
Section~\ref{sec:constructing-knots}.

The next step in studying the algebra of KTGs is to find the relations
in the algebra, and, in particular, the \emph{elementary} relations
(Section~\ref{sec:relat-centr-quest}): pairs of sequences of
operations that produce the same output on all inputs.  Further
justification for calling these relations ``elementary'' comes from
the fact that in other spaces that support the same operations the
elementary relations are automatically satisfied, while other
relations give us non-trivial equations to solve. Tracing out the
track of the KTG as it evolves through a sequence of operations, we
construct a \emph{movie surface} (Section~\ref{sec:movies}), a
decorated simple 2-polyhedron.  If two different sequences of
operations generate the same movie surface, then they are universally
equal.  The converse is not quite true.

In fact, movie surfaces are a special case of \emph{shadow diagrams},
as we briefly discuss in Section~\ref{sec:shadows}.  Shadow diagrams
for 3-manifolds and for 4-manifolds bounded by 3-manifolds were
introduced by Turaev.  They were initially introduced to describe
links inside circle bundles over a surface~\cite{Turaev:ShadowLinks}.
The construction was later generalized to allow descriptions of all
3-manifolds~\cite{Turaev:TopologyShadows, Turaev:quantum}.  From a
quantum topology point of view, shadow diagrams encapsulate the
algebra of quantum 6j-symbols in a concise way; more classically, a
shadow diagram of a 3-manifold can be thought of as an analogue of a
pair of pants decomposition of a surface.

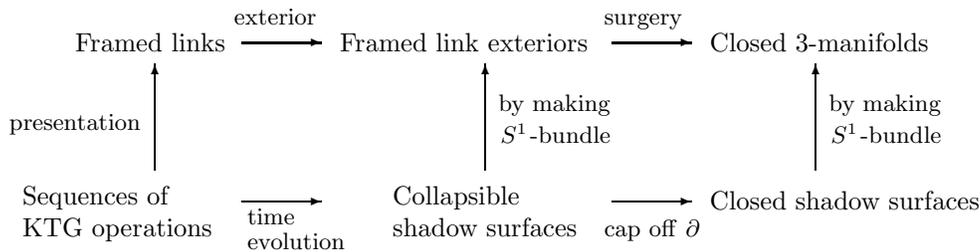
\begin{figure}[htpb]
\begin{center}\vspace{-5mm}
\begin{picture}(320,100)
 \put(0,70){\small Framed links}
 \put(63,73){\vector(1,0){30}}
 \put(60,80){\footnotesize exterior}
 \put(100,70){\small Framed link exteriors}
  \put(203,73){\vector(1,0){30}}
 \put(200,80){\footnotesize surgery}
 \put(240,70){\small Closed 3-manifolds} 
 \put(-25,40){\footnotesize presentation}
 \put(30,25){\vector(0,1){40}}
 \put(160,35){\shortstack[l]{\footnotesize by making \\
                             \footnotesize $S^1$-bundle}}
 \put(155,25){\vector(0,1){40}}
 \put(285,35){\shortstack[l]{\footnotesize by making \\
                             \footnotesize $S^1$-bundle}}
 \put(280,25){\vector(0,1){40}}
 \put(-20,0){\shortstack[l]{\small Sequences of \\
                          \small KTG operations}}
 \put(63,13){\vector(1,0){30}}
 \put(65,-5){\shortstack[l]{\footnotesize time \\ 
                           \footnotesize evolution}}
 \put(120,0){\shortstack[l]{\small Collapsible \\ 
                            \small shadow surfaces}}
 \put(203,13){\vector(1,0){30}}
 \put(200,0){\footnotesize cap off $\partial$}
 \put(240,10){\small Closed shadow surfaces}
\end{picture}
\end{center}
\caption{\label{fig:shadow-summary}%
  A summary of the relations between links, 3-manifolds, KTGs and
  shadow surfaces}
\end{figure}

In summary, every KTG presentation of a knot gives a certain abstract
surface representing that knot complement.  Thus it turns out that our
intrinsically 3-dimensional representation of knot theory can
also be encoded in 2-dimensional terms.  The relation between the
various spaces and constructions is summarized in
Figure~\ref{fig:shadow-summary}.

This paper is an exposition of results that were, at least implicitly,
previously known.  Rather, our aim is to give an exposition of the
relationship between shadow surfaces and knotted trivalent graphs.  In
a future paper, we will prove some theorems that came out of this
work, including the relationship between hyperbolic volume and minimal
complexity shadow diagrams representing a knot.  In addition, the
unifying framework of KTGs is part of an ongoing project to find new
combinatorial presentations of knots, which may be more algebraically
manageable than the full-blown algebra of KTGs.  See
Section~\ref{sec:directions} for more on both of these.

\subsection{Acknowledgements}
A very large part of the credit for this work must go to Dror
Bar-Natan, with whom I have had a long, productive, and fun
collaboration.  It was a great experience coming up with the algebra
of KTGs with him and puzzling over the wide variety of mysterious
relations that we found together.  In this papers,
Sections~\ref{sec:KTG} through~\ref{sec:relat-centr-quest} are joint
work with him.  In addition, I would like to thank Riccardo Benedetti,
Francesco Costantino, Robion Kirby, Tomotada Ohtsuki, Carlo Petronio,
A. Referee, Dale Rolfsen, Chung-chieh Shan, Vladimir Turaev, and
Genevieve Walsh for a great many helpful conversations and comments.  This
work was supported by an NSF Postdoctoral Research Fellowship, a JSPS
Postdoctoral Research Fellowship, and by BSF grant \#1998119.

\section{The space of knotted trivalent graphs}
\label{sec:KTG}

There are relatively few operations on knots.  Traditional operations,
like connect sum, cabling, or general satellites, when applied to
non-trivial knots, never yield hyperbolic knots, but almost all knots
(by any reasonable measure of complexity) are hyperbolic.  Thus
typical knot operations are no good for breaking the vast majority of
knots down into simpler pieces.  To fix this, we will pass to the
larger space of knotted trivalent graphs.  This space will allow more
operations; enough to generate all knots from a few simple generators.

``Graphs'' in this paper might more properly be called ``1-dimensional
complexes''; that is, they may have multiple edges, self loops, and
circle components.

\begin{definition}\label{def:framed-graph}
  A \emph{framed graph} or \emph{fat graph} is a thickening of an
  ordinary graph into a surface (not necessarily oriented): the
  vertices are turned into disks and the edges are turned into bands
  attaching to the disks.  More abstractly, a framed graph is a
  1-dimensional simplicial complex~$\Gamma$ together with an embedding
  $\Gamma\hookrightarrow\Sigma$ of~$\Gamma$ into a surface~$\Sigma$ as
  a spine; more combinatorially, a framed graph is a graph with a
  cyclic ordering on the edges incident to each vertex and $1$ or $-1$
  on each edge (representing a straight or flipped connection,
  respectively), modulo reversing the ordering at a vertex and
  negating the elements on the adjoining edges.
\[
\mathgraph{draws/framed.0}\leadsto\mathgraph{draws/framed.1}
  \quad
\mathgraph{draws/framed.2}\leadsto\mathgraph{draws/framed.3}
  \quad
\mathgraph{draws/framed.10}\leadsto\mathgraph{draws/framed.11}
\]
\end{definition}

The notion of spines comes from PL topology:
a \emph{spine} of a simplicial complex~$Y$ is a
subcomplex~$X$ of $Y$ onto which $Y$ collapses, where collapsing means
successively removing pairs of a $k$-simplex~$\Delta^k$ and a
$(k+1)$-simplex~$\Delta^{k+1}$, where $\Delta^{k+1}$ is the unique $(k+1)$-simplex
having $\Delta^k$ on its boundary.

\begin{figure}[tbp]
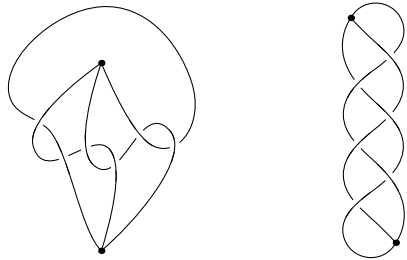

\[
  \mathgraph{draws/examples.0}\qquad\qquad\mathgraph{draws/examples.1}
\]
\caption{\label{fig:examples}%
Two knotted theta graphs, one of which is unknotted.  Which one is it?}
\end{figure}

For instance, the two framed graphs with spine $\Gamma=S^1$ are the annulus
and the M{\"o}bius band.

\begin{definition}
  A \emph{knotted trivalent graph} (KTG) is a trivalent framed
  graph~$\Gamma$ embedded (as a surface) into $\RR^3$, considered up to
  isotopy.
\end{definition}

\begin{definition}
  An \emph{unknotted KTG} is a planar KTG: a KTG which factors through
  an embedding of~$S^2$ in $\RR^3$.
\end{definition}

In diagrams, KTGs will be drawn with the blackboard framing unless
explicitly indicated otherwise.

\section{Elementary operations}
\label{sec:KTGops}

The advantage of KTGs over knots or links is that they support many
operations.  By an \emph{operation} we mean a function from input
trivalent graphs to an output graph, depending only on combinatorial
choices (and no topological data).

\begin{itemize}
\item \emph{Unzip} takes an edge which connects two distinct vertices
  and splits it in two, as though you were unzipping a zipper along
  the edge.  It reduces the number of vertices by~2.
\[
\mathgraph{draws/ops.0} \longmapsto \mathgraph{draws/ops.1}
\]
\item \emph{Bubbling} adds a small loop along an edge of the KTG.
  It increases the number of vertices by~2.
\[
\mathgraph{draws/ops.20} \longmapsto \mathgraph{draws/ops.21}
\]
\item \emph{Connect sum} takes two input trivalent graphs, with a
  chosen edge on each one (and a chosen side of each edge), and
  splices the two edges together.  Note that this operation has two,
  independent, inputs (else the operation depends on a choice of an
  arc and so is not well-defined in the sense above).  Bubbling is
  equivalent to connect sum with an unknotted theta graph~$\thetagraph$.

\[ 
\Bigg\langle\,\mathgraph{draws/ops.10}\;,\;\mathgraph{draws/ops.11}\,\Bigg\rangle
  \longmapsto
\mathgraph{draws/ops.12}
\]  
\item The identity for connect sum is the \emph{unknot}.
\[
\mathgraph{draws/ops.30}
\]

\end{itemize}
Unzip and bubbling operate on one input graph, while connect sum
operates on two inputs and the unknot has no inputs.

\begin{exercise} Check that these operations are well-defined.  In
  particular, connect sum is well-defined for the same reason that
  the connect sum of two knots is well-defined.  (Hint: you can
  shrink one of the trivalent graphs into a small ball and slide it
  along the other.)
\end{exercise}

Note that ``operations'' on knot diagrams such as changing a crossing
or self connect sum are not operations in our sense, since they depend
not only on combinatorial data, but also on a topological choice of
\emph{how} to perform the move.

Also, we only consider these operations in the forward direction.  
The reverse operations cannot always be performed, and so are not KTG
operations, in that they do not act on the space of all KTGs with a
given underlying trivalent graph.

These operations do not suffice to construct interesting knots; for
instance, every KTG constructed with these operations will be planar
(unknotted).  In addition to these operations, we will use three
elementary KTGs as generators:

\begin{itemize}
\item The unknotted \emph{tetrahedron}.
\[
\mathgraph{draws/ops.50}
\]
\item The two minimally twisted \emph{M{\"o}bius bands}, with a positive
  (resp. negative) half twist.  There is no blackboard framing for these
  M{\"o}bius bands, so the half twist is indicated by a little local
  picture of the surface.
\[
\mathgraph{draws/ops.60}\qquad\mathgraph{draws/ops.61}
\]
\end{itemize}

\begin{remark}
  The distinction between operations on KTGs and generators for the
  algebra of KTGs is somewhat a matter of taste.  For instance, the
  unknot can be constructed by unzipping the tetrahedron twice, so
  need not be included in the algebra; however, just as the unit of a
  ring is generally considered to be part of the structure of a ring,
  it is more natural to consider the unknot as a part of the algebraic
  structure of KTGs.
\end{remark}

\section{Constructing knots}
\label{sec:constructing-knots}

The KTG operations we have defined are quite powerful.  In particular,
they generate knot theory.

\begin{theorem}\label{thm:generate}
  Any KTG can be obtained from the unknotted (planar) tetrahedron and
  the two minimally twisted M{\"o}bius bands using unzip, connect sum,
  bubbling, and the unknot.
\end{theorem}

Each of the presentations of knots listed at the beginning of the
paper can be turned into a sequence of KTG operations and thus gives a
proof of Theorem~\ref{thm:generate}.  In the proof of the theorem
below we will show how to turn a knot diagram into a sequence of KTG
operations.  But first, let us define some \emph{composite operations}
formed by composing the elementary operations.
\begin{itemize}
\item A connect sum along an edge followed by an unzip along one of
  the two newly created edges gives a \emph{vertex connect sum}: take
  a vertex in each of two KTGs and join the incoming edges pairwise.
  For a given matching of the edges, there are several ways of
  performing this operation as a sequence of elementary operations,
  all yielding identical results.
\[
\Biggl\langle\mathgraphsmall{draws/ops.40}\;,\;\mathgraphsmall{draws/ops.41}\Biggr\rangle
  \mapstochar\xrightarrow{\text{connect sum}}
\mathgraphsmall{draws/ops.42}
  \mapstochar\xrightarrow{\text{unzip}}
\mathgraphsmall{draws/ops.43}
\]
\item More generally, there is a notion of \emph{tree connect sum}:
  for any KTGs $K_1, K_2$, and isomorphic open subsets $T_1$, $T_2$ of
  the skeletons of $K_1$, $K_2$ which are homeomorphic to trees, we
  can join $T_1$ and $T_2$ by a connect sum and repeated unzips.
  Topologically, you can think of straightening out the two subtrees
  into some standard embedding in a ball (which you can always do,
  since any embedding of a tree in $\RR^3$ is trivial) and then
  removing the balls from both $K_1$ and $K_2$ and identifying the two
  sphere boundary components, connecting corresponding points on the
  boundaries.
\[
\Biggl\langle\mathgraphsmall{draws/ops.70}\;,\;\mathgraphsmall{draws/ops.71}\Biggr\rangle
\longmapsto
\mathgraph{draws/ops.72}
\]
\item Connected sum with M{\"o}bius bands can \emph{change the
    framing} on an edge by any desired (integral or half integral)
  amount.
\item Tree connect sum with the tetrahedron
  \smash[b]{$\mathgraphsmall{draws/misc.40}$} along the dashed subtree
  performs the \emph{Whitehead move}.
  \[
\mathgraph{draws/misc.41} \longmapsto \mathgraph{draws/misc.42}
\]
\end{itemize}

To prove Theorem~\ref{thm:generate}, we use a sweep-out argument based
on a knot diagram.  Mark a circle around a small piece of a knot
diagram and consider the KTG formed by the circle and the portion of
the knot diagram on the inside.  If only one feature of the knot
diagram is enclosed, there are only a few possibilities for the
resulting KTG and it is straightforward to check that they may all be
generated with our generators and operations.  Then sweep this marked
circle outwards by stages, performing a few elementary KTG operations
at each stage.

\begin{figure}
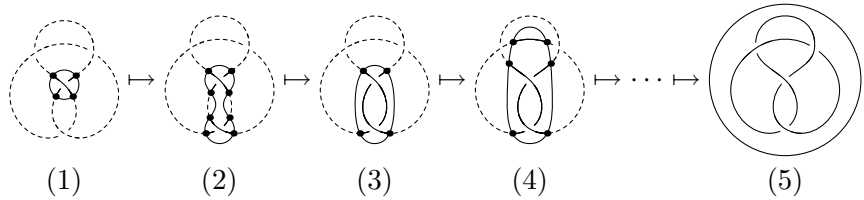

\[
\begin{array}{@{}*8{c@{\;}}c@{}}
\mathgraphsmall{draws/figure8.1}&
\mapsto&
\mathgraphsmall{draws/figure8.2}&
\mapsto&
\mathgraphsmall{draws/figure8.4}&
\mapsto&
\mathgraphsmall{draws/figure8.5}&
\mapsto\cdots\mapsto&
\mathgraphsmall{draws/figure8.8}\\\noalign{\vskip\jot}
(1)&&
(2)&&
(3)&&
(4)&&
(5)
\end{array}
\]
\caption{\label{fig:figure8}%
Generating the figure 8 knot with KTG operations by sweeping
  a circle outwards.}
\end{figure}

Rather than giving an exhaustive proof, we will work through the steps
in a sweepout of a figure 8 knot in Figure~\ref{fig:figure8}, which
provide a good sampling of the cases.  At each stage, the marked
circle and the portion of the knot inside it are drawn solid.
\begin{enumerate}
\item By changing framings, we may turn the unknotted tetrahedron into
  the crossed tetrahedron $\crosstet$.  This is what we obtain by
  drawing our initial marked circle around a crossing, as on the left
  of Figure~\ref{fig:figure8}.
\item At the next step, we push the circle across another crossing.
  We can achieve this by taking connect sum with a second crossed
  tetrahedron.
\item Pushing the circle across a maximum of a strand (looking outwards
  from the strand) is achieved by an unzip move.  Two of these happen
  in the next step.
\item Pushing the circle across a minimum is a bubbling move.
\item We proceed in this way, using one crossed tetrahedron per
  crossing of our diagram, until we have pushed the circle all the way
  out.
\end{enumerate}

\parpic[r]{$\mathgraphsmall{draws/figure8.12}$}
We have succeeded in generating our knot, plus a trivial unknot
component.  There are a few ways to avoid this extra component.  For
instance, in the intermediate steps, one might simply not include the
portions of the circle that touch the exterior region.
This version of step (2) in Figure~\ref{fig:figure8} shown at the right.

Alternatively, we could modify the very last unzip move before
we push the marked circle completely off the knot diagram.
\[
\mathgraphsmall{draws/figure8.7} \,\longmapsto\,
\mathgraphsmall{draws/figure8.13} \,\simeq\,
\mathgraphsmall{draws/figure8.0}
\]

With a similar procedure, we can generate any KTG.  We end up using at
most one tetrahedron per trivalent vertex and one tetrahedron per
crossing in our diagram.

\begin{exercise}\label{xrcs:crossed-tet}
  Show that the crossed tetrahedron is equal to an unknotted tetrahedron
  with the following framing changes:
\[
\raisebox{.5ex-30bp}{\hbox{\graph{draws/misc.15}}} =
  \raisebox{.5ex-30bp}{\hbox{\graph{draws/misc.16}}}
\]
  Hint: start by showing the more elementary equality:
\[
\mathgraphsmall{draws/misc.21} = \mathgraphsmall{draws/misc.20}
\]
\end{exercise}

\section{Relations and elementary relations}
\label{sec:relat-centr-quest}

\begin{figure}
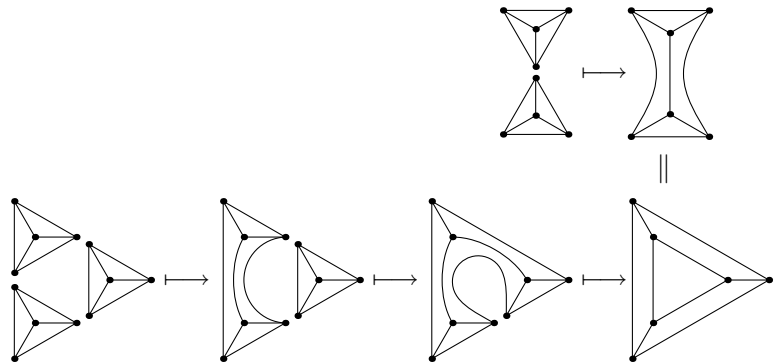

\[
\begin{split}
\mathgraphsmall{draws/rels.0}\longmapsto&\mathgraphsmall{draws/rels.1}\\
&\quad\|\\
\mathgraphsmall{draws/rels.10}\longmapsto
\mathgraphsmall{draws/rels.11}\longmapsto
\mathgraphsmall{draws/rels.12}\longmapsto
&\mathgraphsmall{draws/rels.13}
\end{split}
\]
\caption{The pentagon relation: an unknotted triangular prism can be
  made in two ways.  It may be obtained by the vertex connect sum of
  two tetrahedra (top); alternatively, it may be obtained
  by the vertex connect sum of three tetrahedra followed by an unzip
  (bottom).\label{fig:pentagon}}
\end{figure}

\begin{figure}
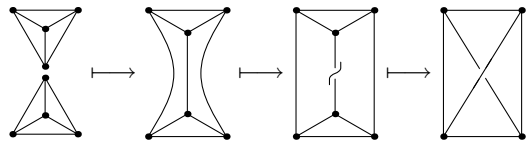

\[
\mathgraphsmall{draws/rels.0}\longmapsto
\mathgraphsmall{draws/rels.1}\longmapsto
\mathgraphsmall{draws/rels.20}\longmapsto
\mathgraphsmall{draws/rels.21}
\]
\caption{The hexagon relation: unzipping a triangular prism (made from
  two tetrahedra) in a twisted fashion yields a twisted tetrahedron.%
\label{fig:hexagon}}
\end{figure}

Knot theory is not freely generated by the tetrahedron and M{\"o}bius
strip; there are some relations.  Two of them are the pentagon
(Figure~\ref{fig:pentagon}) and the hexagon (Figure~\ref{fig:hexagon}).

The names come from the ``pentagon'' and ``hexagon'' relations in the
theory of non-associative tangles~\cite{BarNatan:NAT}, which become
the above relations when non-associative tangles are interpreted as
sequences of KTG operations in a natural way.  In non-associative
tangles, the pentagon and hexagon are the essential identities, the
ones that require significant work to solve; there are, however, many
more identities which are ``elementary'', in that they are
automatically satisfied in the framework.  The same thing is true in
the theory of knotted trivalent graphs.  This leads us to the central
question of this paper, which we can phrase in several different ways.

\begin{question}\label{quest:central}
  What are the elementary relations in the algebra of KTGs?  Which
  composite KTG operations produce the same output for all inputs?
  Which KTG operations are associative?  What relations in the algebra
  of KTGs are automatically true in other natural spaces supporting
  the trivalent graph operations (unzip, bubbling, connect sum,
  unknot)?  What is the operad%
  \footnote{To be precise, this is be a typed operad, whose types are
    the framed trivalent graphs.}
  of KTG operations?
\end{question}

Here are some examples of elementary equalities.  Each elementary
equality equates two composite KTG operations that take the same
number and types of KTGs as input.
\begin{itemize}
\item Any two operations performed on disjoint pieces of a knotted
  graph commute with each other.
\item The connect sum operation on knots is commutative and associative.
\item The parallel operation on links (replacing a knot component with
  a 2 components parallel with respect to the framing) is
  cocommutative and coassociative.
\item The operations of vertex connect sum and tree connect sum are
  well-defined: they do not depend on which of the matching edges we
  connect sum and the order in which we perform the unzips.
\item The left hand (``3'') side of the pentagon equation above does
  not depend on the choices.  We took 3 tetrahedra and did 2 vertex
  connect sums followed by an unzip.  The resulting triangular prism
  has a symmetry that the construction did not have.  This is not an
  accident.  In fact, if we take \emph{any} three knotted tetrahedra
  (with an identification of the underlying graphs with the underlying
  graph of the three standard tetrahedra above) and do the same
  sequence of operations (vertex connect sums between two pairs
  followed by an unzip) in the three possible ways, the result depends
  only on the original knotting of the tetrahedra and not on the order
  of operations.
\item The algorithm for turning a knot diagram into a sequence of KTG
  operations depended on a choice of how to sweep out the knot diagram
  by circles.  In fact, the result is independent of the sweep out,
  even if the crossed tetrahedra $\crosstet$ are replaced by \emph{arbitrary}
  knotted tetrahedra.
\end{itemize}

Some of these identities towards the end of the list may surprise you.
They can all be given good topological explanations.  For instance,
the description of tree connect sum in terms of gluing balls is
clearly independent of choices.  The case of the ``3'' side of the
pentagon can be reduced to some tree connect sums using the following
exercise.

\parpic[r]{$\tripletheta$}
\begin{exercise}
  Express the combination of three tetrahedra as on the ``3'' side of
  the pentagon equation as a tree-connect sum with the graph at right
  along the dashed subtrees: one subtree is connected to each of the
  input tetrahedra.  In other terms, we can straighten out the
  subtrees being acted on in each of the three input tetrahedra before
  performing the operation.  The symmetry of the operation then
  becomes obvious.
\end{exercise}

\parpic[r]{\includegraphics[width=2cm]{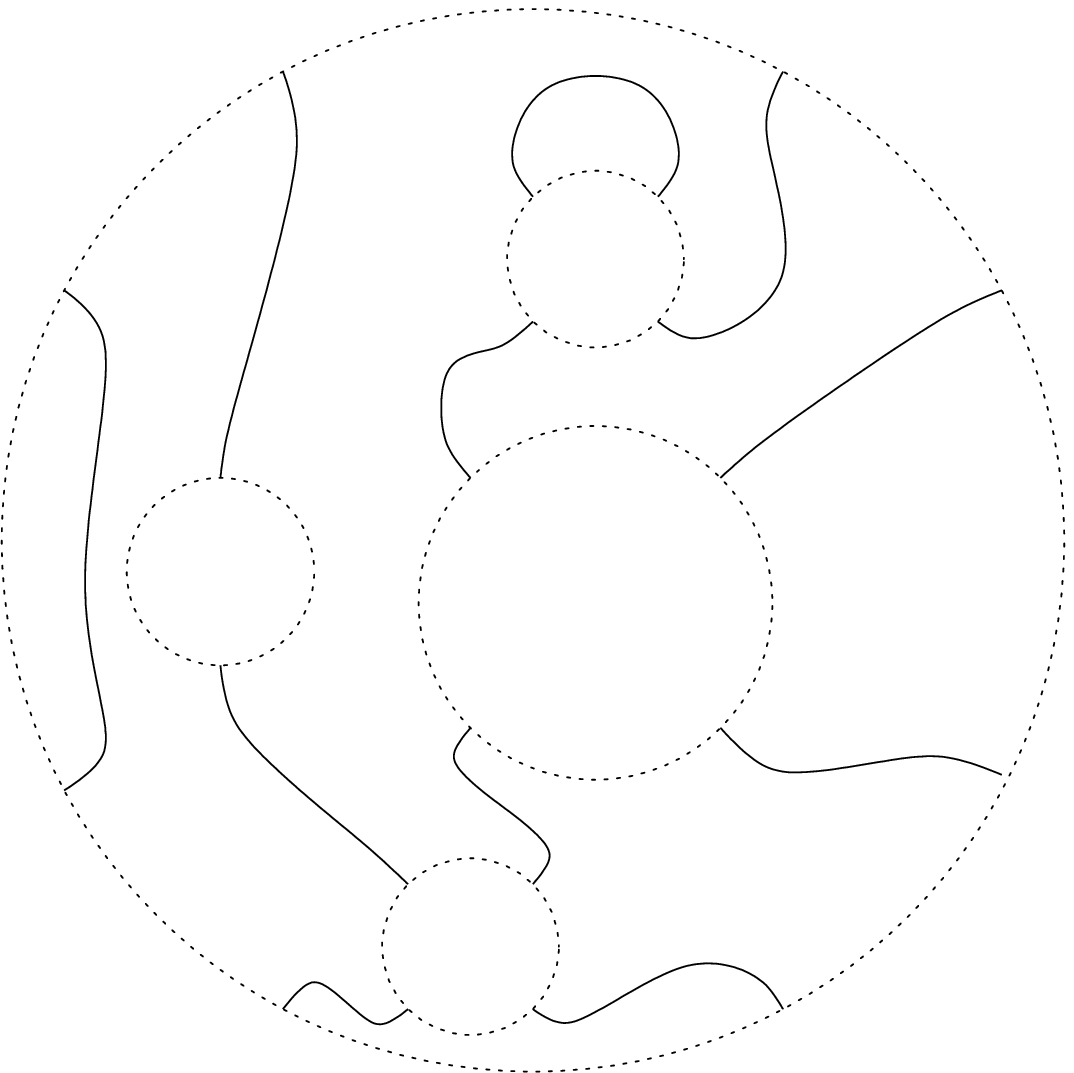}} The
``knot diagram'' case is more involved.  In the case that the outer
circle of the small tetrahedron is unknotted this operation amounts to
gluing together tangles with 4 legs following the pattern of the knot
diagram; however, the operation is well-defined in general.  The
general statement is that knotted trivalent graphs with a
distinguished cycle form a planar algebra in a sense similar to that
of Vaughan Jones~\cite{Jones:PlanarAlgebraI}.  (This structure has
also been called a ``spider''~\cite{Kuperberg:96spiders} or a
``spherical category''~\cite{BarrettWestbury:99spherical}.)  Briefly,
for every arc diagram as on the right, there is a well-defined
operation that takes as input KTGs corresponding to the interior
cycles (with a distinguished cycle containing a corresponding set of
vertices) and produces as output a KTG with a distinguished cycle.
The proof would take us too far afield, but it follows from more
general statements below.

\section{Movie surfaces}\label{sec:movies}

To give some order to the zoo of seemingly ``elementary'' equivalences
between composite knotted trivalent graph operations, we add an extra
dimension.  Consider making a movie of the graph as it evolves in
$S^3$.  The graph traces out a surface with some simple
singularities in a 4-dimensional space.  We can continue this surface
across each of the elementary operations and generators to form a
continuous singular surface, as illustrated in
Table~\ref{tbl:surfaces}.  For instance, to do an unzip move along an
edge, first shrink the edge to zero length and then split the two new
branches apart.  The resulting movie is shown in the first line of the
table.\eject

\begin{table}[tbp]
\begin{center}
\begin{tabular}{lcc}
\toprule
Name & Operation & Movie surface \\
\midrule

Unzip &
$\mathgraph{draws/ops.0}\mapsto\mathgraph{draws/ops.1}$ &
$\mathgraph{draws/ops3d.0}$ \\

Bubbling &
$\mathgraph{draws/ops.20} \mapsto \mathgraph{draws/ops.21}$ &
$\mathgraph{draws/ops3d.20}$ \\

Connect sum &
$\mathgraph{draws/ops.10},\mathgraph{draws/ops.11} \mapsto
\mathgraph{draws/ops.12}$ &
$\mathgraph{draws/ops3d.10}$  \\

Unknot &
$\mathgraph{draws/ops.30}$ &
$\mathgraph{draws/ops3d.30}$ \\

M{\"o}bius band &
$\mathgraph{draws/ops.60}$ &
$\mathgraph{draws/ops3d.60}$ \\

Tetrahedron &
$\mathgraph{draws/ops.50}$ &
$\mathgraph{draws/ops3d.50}$ \\
\bottomrule
\end{tabular}
\end{center}
\caption{The surface created from a movie of the knotted trivalent
  graph operations and generators.  The surfaces depicted all lie in a
  local 3-dimensional slice of a 4-dimensional space; although the
  pictures appear 3-dimensional, they are really pictures of surfaces
  in 4 dimensions.\label{tbl:surfaces}}
\end{table}

Some comments on this table:\nopagebreak[4]
\begin{itemize}
\item These are surfaces in a 4-dimensional space (the evolving
  $S^3$), although the depicted portions lie inside a 3-dimensional
  slice.  This is related to the fact that our knotted graph
  operations locally lie in a plane, as does the unknotted
  tetrahedron (but not the M{\"o}bius strip; see below).
\item The unzip and bubbling surfaces lie in $S^3\times I$, with one
  input $S^3$ (at the bottom of the picture) and one output $S^3$ (at
  the top).  Outside of the portion depicted, there is no change to
  the graph, so the corresponding surface is a product.
\item The connect sum operation has two input $S^3$'s, depicted as
  two spheres in the picture.  (The interiors of these spheres are not
  included in the ambient 4-manifold.)  As before, the output $S^3$ is
  at the top of the picture.  This surface lies inside a
  thrice-perforated $S^4$.
\item The unknot, M{\"o}bius strip, and unknotted tetrahedron have no
  inputs.  The corresponding surfaces lie in $B^4$.
\item The M{\"o}bius strip is the only one of these operations or
  generators that is not essentially planar.  The underlying unframed
  knotted graph of the M{\"o}bius strip is a simple loop, so it bounds a
  disk in $B^4$, just like the unknot.  To distinguish it from the
  unknot we must somehow record the framing.  This will be explained
  more in Section~\ref{sec:shadows}, but for now we will just attach
  the framing information to the surface itself, as indicated in the
  diagram.  (The convention in diagrams is that this framing
  information, the ``gleams'', may be attached to the surface at
  multiple points; each surface component has a total gleam, which is
  the sum of the gleams attached at various points.  In particular,
  if no gleams are indicated on a surface component the total gleam
  is~0.)
\end{itemize}

\begin{figure}
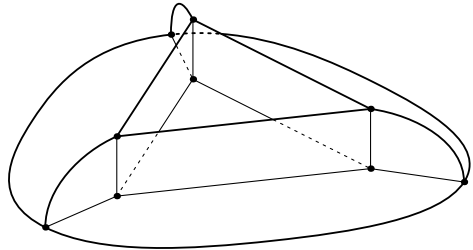

\[
\mathgraph{draws/three.1}
\]
\caption{The movie surface for the ``3'' side of the pentagon.
  The picture lies in $B^4$, with the edges drawn with thick lines
  lying on the boundary of $B^4$.\label{fig:three-1}}
\end{figure}

For a somewhat more involved example of the surface constructed by the
moving knotted trivalent graph, Figure~\ref{fig:three-1} shows the
surface associated to the ``3'' side of the pentagon.  Observe that
the surface has the symmetries of the triangular prism; that is, it
has the symmetries that were missing from the previous description
with KTG operations.

\begin{exercise}
  Convince yourself that Figure~\ref{fig:three-1} is the movie surface
  representing the sequence of operations on the bottom line of
  Figure~\ref{fig:pentagon}.  The action in the movie is from
  the inside out: start with small neighborhoods of each vertex
  (corresponding to 3 unknotted tetrahedra) and push the neighborhoods
  outwards until they meet in pairs (performing the vertex connect sum
  moves) and then encompass the remaining edge (performing the unzip
  move).
\end{exercise}

\begin{figure}
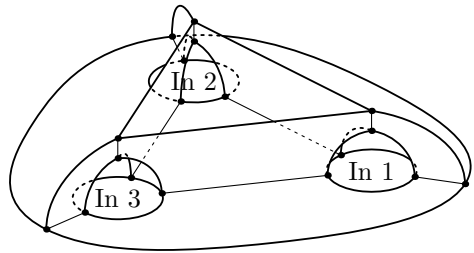

\[
\mathgraph{draws/three.2}
\]
\caption{The movie surface
  for the operation which takes three knotted tetrahedra and attaches
  them in the same pattern as on the ``3'' side of pentagon.%
\label{fig:three-2}}
\end{figure}

The movie surface for the corresponding operation on 3 knotted
tetrahedra (which takes the tetrahedra and connects them according to
the same pattern) is shown in Figure~\ref{fig:three-2}.  This surface
lives inside $B^4$ with three balls removed.  The three input
components of the boundary are labelled ``In 1'', ``In 2'', and ``In
3''; each of these input components lie on one of the removed balls.
The picture suggests that the input tetrahedra are unknotted, but in
fact this surface can be embedded (uniquely!) in $B^4$ minus 3 balls
so as to match any 3 given input knotted tetrahedra.

\begin{figure}
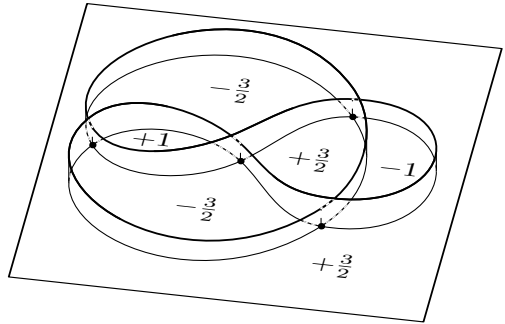

  \[
  \mathgraph{draws/fig8_3d.0}
  \]
  \caption{The movie surface corresponding to the sequence of KTG
    operations used to construct the figure 8 knot in
    Figure~\ref{fig:figure8}.%
  \label{fig:figure8_3d}}
\end{figure}

More generally, consider the movie surface obtained from the sequence
of KTG operations used to construct a knot (plus an extra unknotted
component) from a knot diagram in
Section~\ref{sec:constructing-knots}.  Tracing out the surface through
the sequence of KTG operations, we find that the movie surface is
constructed by taking a disk and attaching an annulus along the curve
of the knot diagram, with the gleams
\[
\mathgraph{draws/misc.17}
\]
on the regions of the knot diagram around each vertex.  Each
complementary region to the knot diagram gets assigned several gleams,
one for each incident vertex; add up the gleams for each vertex.
(These gleams come from Exercise~\ref{xrcs:crossed-tet}).  If we
forget the gleams, the surface is the mapping cylinder of the map from
$S^1$ to the disk specified by the knot diagram.  The resulting
surface coincides with Turaev's ``shadow cylinder'' on the knot.  See
Figure~\ref{fig:figure8_3d} for the example of the figure 8 knot.
Note that this is a picture of abstract surface which does not embed
in $\RR^3$; the transverse intersections of surfaces which appear in
the diagram are only artifacts of the immersion.  There are two
boundary components of this surface (including the external square),
corresponding to the fact that we constructed a 2-component link (with
one unknotted component) in Figure~\ref{fig:figure8}.

\begin{observation}
  When we pass from a sequence of KTG operations to the associated
  surface, we often see additional symmetries.  We do not get the same
  surface for two topologically different sequences of KTG operations.
\end{observation}

Thus the movie surface seems to be a good representation of a sequence
of KTG operations.  This observation will be justified below,
where we will see that a KTG may be reconstructed from its movie
surface.

\begin{exercise}
  Check that each of the sequences of equivalent KTG operations listed
  in Section~\ref{sec:relat-centr-quest} yield homeomorphic surfaces,
  with one exception: coassociativity of the parallel operation.
\end{exercise}

\section{Shadow Surfaces}\label{sec:shadows}

The movie surfaces constructed in the previous section are, in fact,
special cases of \emph{shadow surfaces}, which we will now introduce.
The key observation to make about each of the drawings in
Table~\ref{tbl:surfaces} is that the ambient 4-manifold collapses onto
the union of the movie surface and the $S^3$ input(s).
To see this collapsibility, first note that outside a neighborhood of
the operation, the movie surface is constant in time, and so the
4-manifold away from the movie surface may be collapsed onto the
input.  Furthermore, the operations all take place inside a
3-dimensional slice; the 4-manifold outside of the 3-dimensional slice
can again be collapsed.  So we are left with showing that, for each
surface~$\Sigma$ in Table~\ref{tbl:surfaces} considered as a
3-dimensional surface, $\RR^3$ collapses down onto~$\Sigma$.  But the
complementary regions of each~$\Sigma$ are all 3-balls touching the
output surface in just one connected 2-balls, so we can collapse each
3-ball, as desired.

Note that the argument above does not apply if, e.g., we
turned the unzip surface upside down (one complementary region
touches the output in two places), turned the bubble surface upside
down (one complementary region does not touch the output), or
considered a self-connect sum with only one input (one complementary
region is not a 3-cell); corresponding to the fact that none of these
operations are topologically well-defined inside $S^3$.  (Shadow
surfaces do give a topological interpretation to each of these
surfaces, but the result is a KTG in some 3-manifold, not generally
$S^3$.)

A somewhat stronger statement is that in each case the ambient
4-manifold is diffeomorphic to a collar over the $S^3$ input(s),
together with a regular neighborhood of~$\Sigma$ in a 4-manifold, with
the surface embedded in a locally flat way. (Without this caveat, the
surface might locally be the cone over some non-trivial knot in
$S^3$.)

The situation is somewhat easier to think about in the case when there
are no inputs, so we have a KTG presentation of a knotted
graph~$\Gamma$.  In this case, all of $B^4$ collapses down onto the
movie surface~$\Sigma$ and, to reconstruct the pair $(B^4, \Gamma)$
from~$\Sigma$, we just need to describe how to ``thicken'' $\Sigma$
into its regular neighborhood.  This thickening should take a surface
(with singularities) and yield a smooth 4-manifold.  Away from the
singularities, we take a disk bundle over the surface, and continue in
a natural way over the singularities.  If we are only interested in
3-manifold topology, we may consider the circle bundle at the boundary
of the disk bundle.

To understand this situation a little better, let us consider the
situation one dimension down: instead of constructing 3-manifolds from
singular surfaces, let us construct 2-manifolds from graphs.  (See
Figure~\ref{fig:shadow-summary} or Table~\ref{tab:codim-analogy} for a
guide to the analogy.)

\begin{table}[htbp]
\begin{center}
\begin{tabular}{rll}
\toprule
  &
    Dimension 1 (graphs) &
      Dimension 2 (simple surfaces) \\
\midrule
  codim 1 &
    Fat graph &
      Simple spine of 3-manifold\\
  codim 2 &
    3-dimensional handlebody &
      Shadow of slim 4-manifold \\
  $\partial$ of codim 2  &
    Pair-of-pants decomposition &
      Shadow of 3-manifold \\
\bottomrule
\end{tabular}
\end{center}
\caption{An analogy between shadows of 3-manifolds and pair-of-pants
  decompositions of surfaces: the regular neighborhoods of polyhedra
  of varying dimension and codimension.\label{tab:codim-analogy}}
\end{table}

We have already seen one way to construct 2-manifolds from graphs, in
Definition~\ref{def:framed-graph}: The construction of fat (or framed)
graphs from a trivalent graph~$\Gamma$.  This construction is nearly
equivalent to triangulations of closed surfaces: if we glue a disk
onto each circle component of the boundary of a fat graph, we get a
closed surface~$\Sigma$ with $\Gamma$ drawn on it.  Assuming that $\Gamma$ has no
circle components, the dual cell division to~$\Gamma$ is a triangulation
of~$\Sigma$, in the weak sense (for instance, two sides of a triangle may
be glued to each other).

But this construction is in codimension~1 (a graph is thickened into a
surface), while the construction we are interested in is in
codimension~2 (a surface is thickened into a 4-manifold).  For the
proper analogue, we should consider a codimension~2 thickening of a
graph: that is, thickening a graph into a 3-manifold.
A regular neighborhood of a graph~$\Gamma$ in a 3-manifold
is a handlebody; its boundary is a closed surface~$\Sigma$.  If we cut~$\Sigma$
along circles surrounding each edge of~$\Gamma$, we end up with a
thrice-perforated sphere for each vertex of~$\Gamma$.  This is also known
as a pair-of-pants decomposition of~$\Sigma$.  (We are again assuming that
$\Gamma$ has no circle components.)

Moving up one dimension, let us consider the thickening of surfaces.
First let us be precise about the type of singular surfaces we
consider.

\begin{definition}
  A \emph{simple polyhedron} or \emph{simple surface} is a
  2-dimensional simplicial complex in which the link of every point is
  locally homeomorphic to a neighborhood of some point in the (closed)
  cone over the 1-skeleton of a tetrahedron.  The \emph{boundary} of a
  simple surface is the subset of the surface locally homeomorphic to
  a point on the 1-skeleton of the tetrahedron itself.
\end{definition}

The possible local models for the interior of a simple surface are
$\RR^2$, the three page book, and the open cone over a tetrahedron.
The possible local models for the boundary are the upper half plane
and the upper half of the three page book.  All the surfaces
considered in this paper are simple polyhedra.

In codimension 1, if a manifold $M^3$ is a regular neighborhood of a
simple surface $\Sigma^2$, then $\Sigma$ is called a \emph{spine} of $M^3$.
Every 3-manifold with non-empty boundary has a spine.  Spines of
3-manifolds, like fat graphs, are closely related to triangulations:
the dual cell complex to a triangulation of a closed 3-manifold $M^3$
is a spine for $M^3$ minus a neighborhood of each vertex of the
triangulation.  The principal advantage of spines over triangulations
is that it is easier to consider various kinds of degeneracy: the dual
to a triangulation is a \emph{standard spine}, in which all of the
edges of the spine are intervals connecting two vertices and
all of the faces are disks.  In general, a \emph{standard surface} is
a simple surface satisfying this extra condition that the edges be
intervals and the faces be disks.  It turns out that the thickening of
a standard surface into an orientable 3-manifold is completely
determined by the surface, if it
exists~\cite{Casler:ImbeddingTheorem}.

In codimension 2, we may consider a regular neighborhood~$W^4$ of a
2\hyp{}polyhedron $\Sigma^2$ in an orientable 4-manifold.  For reasonable
topological control, we require that~$\Sigma$ be locally flat at
generic points of $\Sigma$, so that $W$ is generically a
disk bundle over~$\Sigma$.  Suppose for the moment that $\Sigma$ is a
smoothly embedded closed ordinary surface.  Then the disk bundle is
determined by its Euler number.%
\footnote{$\Sigma$ itself need not be oriented: a disk bundle over a
  surface has a well-defined Euler number if the total space is
  orientable.}
This number is what we called the ``gleam'' in
Section~\ref{sec:movies} in the discussion of the M{\"o}bius strip.

\parpic[r]{\graph{draws/misc3d.0}}
In general, the normal bundle to each smooth component~$Y$ of~$\Sigma$
is a 2-plane bundle on $Y$.  In order to assign a well-defined Euler
number to this bundle, we need to define a trivialization on the
boundary.  At generic points on the boundary of $Y$, there are two
other regions of~$\Sigma$ attached to the 1-skeleton.  At vertices on
the boundary of $Y$, there is also an opposite region attached, but it
may be ignored: around each boundary component of $Y$, there is a
unique way to map an annulus or a M\"obius strip~$S$ (the hashed
region in the diagram to the right) in to $\Sigma$ so that the
boundary of $Y$ is a spine of~$S$, and $Y$ is otherwise disjoint
from~$S$.

After homotopy, we can assume that $S$ is normal to $Y$; it
then provides two vectors in the normal bundle to $Y$ over each point
in $\partial Y$.  In the case that $S$ is an annulus, we can pick one
of the two vectors and get a trivialization of the normal bundle over
the boundary, and so we can define a relative Euler number relative
to~$\Sigma$.  In case $S$ is a M\"obius strip, we get a
half-trivialization lying between two honest trivializations (each
obtained by inserting a half twist in $S$); we correspondingly define
the relative Euler number to be a half-integer.  It turns out that
these Euler numbers, the gleams, are enough to reconstruct~$W$
from~$\Sigma$~\cite{Turaev:TopologyShadows, Turaev:quantum}.  (Recall
that $W$ is a regular neighborhood of~$\Sigma$.)
Let us state this precisely.

\begin{definition}
  A \emph{shadow surface} is a simple polyhedron (possibly with
  boundary) together with a gleam in $\ZZ$ (resp. $\ZZ + \frac{1}{2}$)
  on each face which has an even (resp. odd) number of M{\"o}bius strip
  components around the boundary.
\end{definition}

\begin{definition}
  A (strict) \emph{shadow} of a 4-manifold~$W^4$ is a closed shadow
  surface~$\Sigma^2$ which is a spine of~$W$, with the gleam on each
  face~$Y$ of~$\Sigma$ equal to the Euler class of the normal bundle
  to~$Y$.
  
  A \emph{(relative) shadow} of a pair~$(W^4,\Gamma)$ of a 4-manifold~$W^4$
  with a (framed) trivalent graph~$\Gamma$ embedded in $\partial W$ is a
  simple polyhedron~$\Sigma^2$ properly embedded in $W$ so that $\partial \Sigma = \Gamma$
  and $W$ collapses onto $\Sigma$.  Recall that the framing really means
  the thickening of $\Gamma$ into a surface; attach $W$ to the original
  graph running along the center of this surface.  The resulting
  slightly larger surface can be used to give gleams to all the
  regions of $W$ in a uniform way.
\end{definition}

\begin{figure}[p]
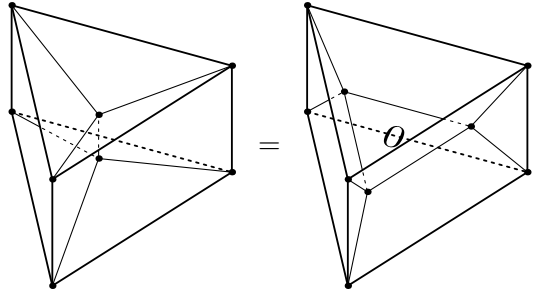

  \[
  \mathgraph{draws/rels3d.0} = \mathgraph{draws/rels3d.1}
  \]
  \caption{The pentagon relation in the shadow world.%
  \label{fig:pentagon-shadow}}
\end{figure}

\begin{figure}
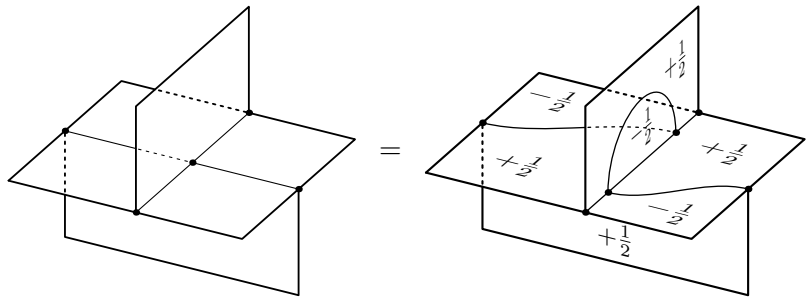

  \[
  \mathgraph{draws/rels3d.10} = \mathgraph{draws/rels3d.11}
  \]
  \caption{The hexagon relation in the shadow world.  On the right
    hand side, the lower rectangle attaches along the path that runs
    on the upper rectangle; the result is not embeddable in $\RR^3$.%
  \label{fig:hexagon-shadow}}
\end{figure}

\begin{remark}
  Our terminology differs somewhat from that of
  Turaev~\cite{Turaev:TopologyShadows, Turaev:quantum}.  Turaev
  defines shadows to be equivalence classes of surfaces modulo some
  moves (including the pentagon and hexagon in
  Figures~\ref{fig:pentagon-shadow} and~\ref{fig:hexagon-shadow}).
  Since these equivalence classes turn out to be nearly the same as
  4-manifolds in the most interesting case, we prefer to reserve the
  name ``shadow'' or ``shadow surface'' for the actual decorated
  surface, rather than the equivalence class.
\end{remark}

\begin{proposition} [Turaev]
  Every shadow surface~$\Sigma$ is the shadow of a 4-manifold~$W$, which is
  unique up to PL homeomorphism preserving $\Sigma$.
\end{proposition}

\begin{remark}
  In the setting of the above proposition, $W$ may admit non-trivial
  self\hyp{}homeomorphisms preserving $\Sigma$.
\end{remark}

\begin{proposition}
  The movie surface associated to a KTG presentation of a knotted
  graph~$\Gamma$ is a collapsible shadow for $(B^4,\Gamma)$.  Furthermore, every
  collapsible shadow surface arises in this way.
\end{proposition}

\begin{proposition}
  Every collapsible shadow surface~$\Sigma$ (with arbitrary gleams) is the
  shadow of $(B^4,\Gamma)$ for some KTG~$\Gamma$.
\end{proposition}

\begin{proof}
  Let $\Sigma$ be a shadow for $(W^4, \Gamma)$.  By hypothesis, $W$
  collapses onto $\Sigma$ and $\Sigma$ collapses to a point, so $W$
  collapses to a point and is therefore $B^4$.
\end{proof}

\begin{remark}
  A KTG presentation is equivalent to its movie surface~$\Sigma$
  together with the horizontal slicing given by time.  This gives a
  particular type of Morse function on $\Sigma$, one which exhibits
  the fact that $\Sigma$ is collapsible.  (Precisely, for $t_0 \leq
  t_1$, $\Sigma_{t_0\leq t\leq t_1}\cup C(\Sigma_{t_0})$ should be
  collapsible, where $C(\Sigma_{t_0})$ is the cone on $\Sigma_{t_0}$.)
\end{remark}

\begin{question}
  Does every collapsible simple 2-polyhedron~$\Sigma$ admit a Morse
  function exhibiting the fact that $\Sigma$ is collapsible?
\end{question}

\begin{definition}
  A \emph{shadow representation} of an oriented  3-manifold~$M^3$ is a
  4-manifold $W^4$ with $\partial W = M$, together with a shadow for $W$.
\end{definition}

Every oriented 3-manifold has a shadow representation: The effect of
attaching a disk to the boundary of a shadow surface is surgery on the
corresponding knot.  Since every link in $S^3$ has a shadow
representation (since it has a KTG presentation) and every 3-manifold
is surgery on a link, we can get shadow representations of any
3-manifold.

On the other hand, not every 4-manifold has a shadow representation.
A shadow representation of a 4-manifold~$W$ provides a handle
decomposition of $W$ with only 0-, 1-, and 2-handles, so we must in
particular have $H^3(W) = H^4(W) = 0$.  This class of 4-manifolds is
an interesting one; for instance, they are the handlebodies that are
homeomorphic to Stein domains~\cite{Gompf:SteinHandlebody}.

\begin{figure}
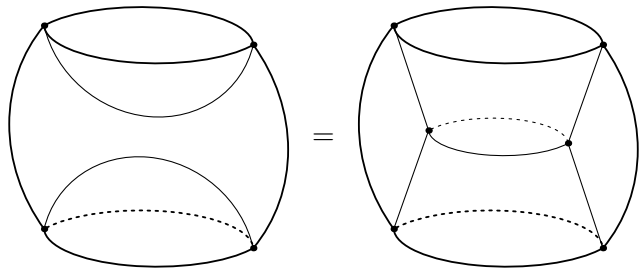

  \[
  \mathgraph{draws/rels3d.30} = \mathgraph{draws/rels3d.31}
  \]
  \caption{The 2-0 move, an additional move needed to be able to
    change the topology of regions of the shadow.%
  \label{fig:2-0-shadow}}
\end{figure}

\begin{question}\label{quest:shadow-calculus}
  Are any two standard shadow surfaces representing the same 4-manifold
  are related by a sequence of pentagon
  (Figure~\ref{fig:pentagon-shadow}), hexagon
  (Figure~\ref{fig:hexagon-shadow}), and 2-0 (Figure~\ref{fig:2-0-shadow}) moves and their inverses?
\end{question}

The answer to Question~\ref{quest:shadow-calculus} is probably ``no''.
There is a well-known conjecture which is closely related to the case
when the shadow is contractible:
\begin{conjecture}[Andrews-Curtis]\label{conj:andrews-curtis}
  Any two contractible simple polyhedra can be
  related by the pentagon, hexagon, and 2-0 moves and their inverses,
  ignoring all gleams.
\end{conjecture}
There is also a generalized Andrews-Curtis conjecture, which drops the
hypothesis that the polyhedra be contractible.  The usual form of the
Andrews-Curtis conjecture is in terms of balanced presentations of the
trivial group, but this version is equivalent.  This conjecture is
generally believed to be false: there are many known potential
counterexamples.  While Conjecture~\ref{conj:andrews-curtis} and
Question~\ref{quest:shadow-calculus} are not directly related (in
either direction), many potential counter-examples to the
Andrews-Curtis conjecture can be turned into potential
counter-examples to Question~\ref{quest:shadow-calculus} above.  For
instance, the examples of Akbulut and
Kirby~\cite{AkbulutKirby:Counterexample} are of this type.
Question~\ref{quest:shadow-calculus} can be thought of as an embedded
version of the Andrews-Curtis conjecture.

\begin{figure}[htbp]
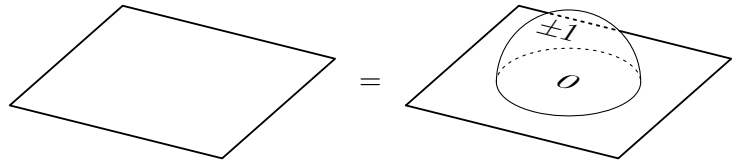

  \[
  \mathgraph{draws/rels3d.20} = \mathgraph{draws/rels3d.21}
  \]

  \caption{$+1$ surgery on an unknot is trivial, as seen in the shadow
  world.%
  \label{fig:surgery-shadow}}
\end{figure}

If we want a shadow calculus for 3-manifolds as opposed to
4-manifolds, we need to allow for some move that changes the
4-manifold.  One such move, which may be used to achieve either the
Kirby I move or the Fenn-Rourke move, is shown in
Figure~\ref{fig:surgery-shadow}.

In joint work with Francesco Costantino, we have shown that this is
sufficient.

\begin{theorem}[Costantino-Thurston, to appear]\label{thm:shadow-calculus}
  Any two simply-connected shadow presentations of the same 3-manifold
  are related by the pentagon, hexagon, 2-0, and $\pm1$-bubble moves.
\end{theorem}

Notice that adding a bubble, and so increasing the
second homology of the surface, is enough to avoid Andrews-Curtis
problems.

\section{What's it good for?}
\label{sec:directions}

Besides providing a unified framework for viewing many different types
of knot representations, the algebra of KTGs as presented in this paper has
many possible applications, to be explored more in future papers.

\subsection{Representations of the algebra of KTGs}
As mentioned above, there are other natural spaces that support the
same set of KTG operations, with the same elementary relations.  We
may call these ``TG-algebras''.  A TG-algebra consists of a space
associated to each (abstract) trivalent graph, and maps between them
corresponding to the elementary operations.  To find a knot invariant
with values in a TG-algebra it suffices to find values for the
tetrahedron and M\"obius strips which satisfy the pentagon and hexagon
relations.

This construction is known for some TG-algebras.  There
is one TG-algebra related to representations of a group~$G$; the
quantum invariants of Reshetikhin and Turaev come from representations
on this space.  For some other spaces, this is apparently new.  For
instance, the space of chord diagrams on the graph, the target
for the universal Vassiliev invariant, supports these same operations;
in these terms, there is an elegant characterization of the Kontsevich
integral.

\begin{proposition}
  The Kontsevich integral $Z$ is the unique universal Vassiliev
  invariant which has an extension to knotted trivalent graphs
  well-behaved under the TG operations, and for which a half-framing
  change acts by multiplication by $\exp(\theta/4)$.
\end{proposition}

Although the knot invariant is unique, the extension to KTGs is not
unique, but the ambiguity is understood.  The condition on the framing
change is just a normalization condition.  This is a
reformulation of previously known results~\cite{drin:quasitri,
  Poirier:ConfigSpace, Lescop:Uniqueness}, but it has extra symmetry.

Another TG-algebra which is a potential target for a KTG invariant is
more-or-less dual to the space of representations mentioned above: the
space associated to a KTG consists of measures on $G^E/G^V$: a copy of
the group~$G$ for each edge, modulo an action for each vertex.  One
reason to be interested in representations on this space is that it
provides a natural setting for trying to understand Witten's
asymptotics conjectures on asymptotic behaviour of quantum invariants.

\subsection{Finding new specializations}
The algebra of knotted trivalent graphs may be too large for some
purposes.  We have only described the operations (up to equivalence)
implicitly, by equivalence of the underlying surfaces.  However, there
are several special cases which are easier to deal with, forming, for
instance, an ordinary algebra with a single associative
multiplication.  Some of these special cases were mentioned at the
beginning of this paper.  The KTG point of view suggests several more
special cases which are worth investigating, including a category of
\emph{annular braids}.

\subsection{Complexity issues}

Every knot diagram with $n$ crossings can be turned into a KTG
presentation with $n$ or fewer tetrahedra.  The converse is not true:
KTG presentations may be much more efficient than knot diagrams.  For
instance, it is easy to construct a sequence of links where the
linking number grows exponentially in the number of KTG operations.

For bounds in the other direction, recall that the Gromov norm of a
knot complement is (up to a constant) equal to the sums of the
hyperbolic volumes of the hyperbolic pieces of the geometric
decomposition of the complement.  (In particular, for an iterated
torus knot the Gromov norm is 0, and for a hyperbolic knot the Gromov
norm is essentially the volume.)  It can be shown that the Gromov norm
is less than a constant times the number of tetrahedra in any KTG
presentation.  So the minimum number of tetrahedra in a KTG
representation of a knot is bounded above by the Gromov norm, and
below by the usual crossing number of the knot.

It is an open question where in between these bounds
the minimum number of tetrahedra lies.  A related result is that any
geometric 3-manifold has a shadow presentation with a number of
vertices at most quadratic in the Gromov
norm.  This is joint work with Francesco Costantino; a paper is in
preparation.

\Addresses\recd

\begin{thebibliography}

\bibitem{AkbulutKirby:Counterexample}
\textbf{Selman Akbulut}, \textbf{Robion Kirby}, \emph{A potential smooth
  counterexample in dimension {$4$} to the {P}oincar\'e conjecture, the
  {S}choenflies conjecture, and the {A}ndrews-{C}urtis conjecture}, Topology 24
  (1985) 375--390

\bibitem{BarNatan:NAT}
\textbf{Dror Bar-Natan}, \emph{Non-associative tangles}, from: ``Geometric
  topology (Athens, GA, 1993)'', Amer. Math. Soc. (1997)  139--183

\bibitem{BarrettWestbury:99spherical}
\textbf{John~W Barrett}, \textbf{Bruce~W Westbury}, \emph{Spherical
  Categories}, Adv. Math. 143 (1999) 357--375

\bibitem{Casler:ImbeddingTheorem}
\textbf{B\,G Casler}, \emph{An imbedding theorem for connected $3$-manifolds
  with boundary}, Proc. Amer. Math. Soc 16 (1965) 559--566

\bibitem{Conway:Enumeration}
\textbf{John~H Conway}, \emph{An enumeration of knots and links, and some of
  their algebraic properties}, from: ``Computational Problems in Abstract
  Algebra (Proc. Conf., Oxford, 1967)'', Pergamon, Oxford (1970)  329--358

\bibitem{drin:quasitri}
\textbf{Vladimir~G Drinfel'd}, \emph{On Quasitriangualar Quasi-{H}opf Algebras
  and a Group Closely Connected with
  {$\mathrm{Gal}(\overline{\mathbb{Q}}/\mathbb{Q})$}}, Leningrad Math. J. 2
  (1991) 829--860

\bibitem{Dynnikov:3Page}
\textbf{Ivan~A Dynnikov}, \emph{Three-page representation of links}, Uspekhi
  Mat. Nauk 53 (1998) 237--238

\bibitem{Gompf:SteinHandlebody}
\textbf{Robert~E Gompf}, \emph{Handlebody construction of {S}tein surfaces},
  Ann. of Math. (2) 148 (1998) 619--693

\bibitem{Jones:PlanarAlgebraI}
\textbf{Vaughan F\,R Jones}, \emph{Planar Algebras, {I}}, Technical report, U.
  C. Berkeley (1999), to appear in New Zealand J. Math., 
  \arxiv{math.QA/9909027}

\bibitem{Kuperberg:96spiders}
\textbf{Greg Kuperberg}, \emph{Spiders for Rank 2 {L}ie Algebras}, Comm. Math.
  Phys. 180 (1996) 109--151, \arxiv{q-alg/9712003}

\bibitem{LM:FramedOriented}
\textbf{Thang T\,Q Le}, \textbf{Jun Murakami}, \emph{The universal
  {V}assiliev-{K}ontsevich invariant for framed oriented links}, Compositio
  Math. 102 (1996) 41--64, \arxiv{hep-th/9401016}

\bibitem{Lescop:Uniqueness}
\textbf{Christine Lescop}, \emph{About the Uniqueness of the {K}ontsevich
  Integral}, Journal of Knot Theory and its Ramifications 11 (2002) 759--780,
  \arxiv{math.GT/0004094}

\bibitem{Poirier:ConfigSpace}
\textbf{Sylvain Poirier}, \emph{The Configuration Space Integral for Links and
  Tangles in {${\mathbb R}^3$}}, Ph.D. thesis, Fourier Institute, University of
  Grenoble (2000), \arxiv{math.GT/0005085v2}

\bibitem{Turaev:TopologyShadows}
\textbf{Vladimir~G Turaev}, \emph{Topology of shadows} (1991), preprint

\bibitem{Turaev:ShadowLinks}
\textbf{Vladimir~G Turaev}, \emph{Shadow links and face models of statistical
  mechanics}, J. Differential Geom. 36 (1992) 35--74

\bibitem{Turaev:quantum}
\textbf{Vladimir~G Turaev}, \emph{Quantum invariants of knots and 3-manifolds},
  W. de Gruyter, Berlin (1994)

\end{thebibliography}
\end{document}